\documentclass[12pt]{amsart}
\usepackage{amssymb}
\usepackage{amscd}
\usepackage{pstricks}
\usepackage{pst-node}
\newtheorem{thm}{Theorem}[section]
\newtheorem{cor}[thm]{Corollary}
\newtheorem{lem}[thm]{Lemma}
\newtheorem{prop}[thm]{Proposition}

\theoremstyle{definition}
\newtheorem{defn}[thm]{Definition}

\newtheorem{ass}{Assumption}

\newtheorem{rem}[thm]{Remark}

\theoremstyle{remark}

\numberwithin{equation}{section}

\newcommand{\norm}[1]{\left\lVert{#1}\right\rVert}

\newcommand{\R}{\mathbb{R}}
\newcommand{\C}{\mathbb{C}}
\newcommand{\Z}{\mathbb{Z}}

\newcommand{\trans}[1]{{}^t\kern-.2em{#1}}
\newcommand{\ytrans}[1]{{}^t\kern-.11em{#1}}
\newcommand{\Trans}[1]{{}^T\kern-.2em{#1}}
\newcommand{\lsup}[2]{{}^{#1}\kern-.1em{#2}}

\newcommand{\Id}{\operatorname{Id}}

\newcommand{\cal}[1]{\mathcal{#1}}

\newcommand{\Ker}{\operatorname{Ker}}
\newcommand{\sign}{\operatorname{sign}}
\newcommand{\Mas}{\mathbf{Mas}}
\newcommand{\M}{\mathbf{\M}}
\newcommand{\Sf}{\mathbf{Sf}}
\newcommand{\Hi}{\mathbb{H}}
\newcommand{\J}{\mathbb{J}}
\newcommand{\graph}{\operatorname{graph}}

\renewcommand{\tilde}[1]{\widetilde{#1}}
\renewcommand{\hat}[1]{\widehat{#1}}

\newcommand{\inner}[2]{\left\langle{#1},{#2}\right\rangle}

\makeatletter
\def\eqnarray{%
   \stepcounter{equation}%
   \def\@currentlabel{\p@equation\theequation}%
   \global\@eqnswtrue
   \m@th
   \global\@eqcnt\z@
   \tabskip\@centering
   \let\\\@eqncr
   $$\everycr{}\halign to\displaywidth\bgroup
       \hskip\@centering$\displaystyle\tabskip\z@skip{##}$\@eqnsel
      &\global\@eqcnt\@ne \hfil$\displaystyle{{}##{}}$\hfil
      &\global\@eqcnt\tw@ $\displaystyle{##}$\hfil\tabskip\@centering
      &\global\@eqcnt\thr@@ \hb@xt@\z@\bgroup\hss##\egroup
        \tabskip\z@skip
      \cr}
\makeatother

\makeatletter
\def\varin{\mathrel{\mathpalette\@varin\relax}}
\def\@varin#1{%
   \hbox{\setbox\z@\hbox{\m@th$#1\cup$}%
       \def\reserved@a{bold}%
       \dimen@\ifx\reserved@a\math@version .3\else .2\fi\p@
       \kern.5\wd\z@\kern-\dimen@
       \vrule\@width2\dimen@\@height1.08\ht\z@\@depth\z@
       \kern-\dimen@\kern-.5\wd\z@
       \box\z@}}
\makeatother

\begin{document}
\title [splitting formula for a spectral flow] 
{Maslov index in the infinite dimension and a splitting formula
       for a  spectral flow}



\author{Kenro Furutani and Nobukazu Otsuki}

\address{Kenro Furutani \endgraf 
Department of Mathematics \endgraf 
Faculty of Science and Technology \endgraf 
Science University of Tokyo \endgraf 
2641 Noda, Chiba 278-8510 Japan. \endgraf}
\email{furutani@ma.noda.sut.ac.jp}
\address{Nobukazu Otsuki \endgraf      
Department of Mathematics \endgraf 
Faculty of Science and Technology \endgraf 
Science University of Tokyo \endgraf 
2641 Noda, Chiba 278-8510 Japan. \endgraf }
\email{nobukazu@ma.noda.sut.ac.jp}

\bigskip
\bigskip
\maketitle
\tableofcontents
\thispagestyle{empty}
\begin{abstract}
  First, we prove a local spectral flow formula (Theorem 3.7)
for a differentiable curve of selfadjoint Fredholm operators.
This formula enables us to prove in a simple way 
a general spectral flow formula
(Theorem 3.8) 
which was already proved in \cite{BF1}.
   Secondly, we prove a splitting formula (Theorem 4.12) 
for the spectral flow of a curve of selfadjoint elliptic operators
on a closed manifold, which we decompose into two parts with commom
boundary.  
Then the formula says that the spectral flow is a sum
of two spectral flows on each part of the separated manifold with
naturally introduced elliptic boundary conditions. 
In the course of proving this formula,
we investigate a property of the Maslov index for paths of Fredholm pairs of 
Lagrangian subspaces.
\end{abstract}

\section{Introduction}

It is well known that the non-trivial component
$\widehat{\mathcal{F}}_*$ of 
the space  $\widehat{\mathcal{F}}=\widehat{\mathcal{F}}_{-}\cup
\widehat{\mathcal{F}}_{*}\cup \widehat{\mathcal{F}}_{+}$ of
bounded selfadjoint Fredholm operators (on a complex or a real Hilbert
space) 
is a classifying space for the $K$-group,
$K^{-1}$ in the complex case and that of $KO^{-7}$ in the real case 
(\cite{AS}).  Also we know, by Bott periodicity theorem, that both of the
fundamental groups of these spaces is isomorphic to $\mathbb{Z}$.
This isomorphism is given by the quantity ``spectral flow'',
although it is not stated explicitly in this paper \cite{AS}.
We have an
``intuitive'' understanding of the spectral flow by saying that the spectral
flow is the difference of the ``net numbers'' of the eigenvalues of the
selfadjoint Fredholm operators which change signs from
minus to plus and from plus to minus, when the parameter of the family
goes from $0$ to $1$.

Since the ``spectral flow'' was treated in the paper \cite{APS}, 
the quantity ``spectral flow'' appears in the various theories 
where it plays important roles, for example, 
in the study of spectral analysis of Dirac operators, 
in the theory of Floer
homology for the study of low dimensional manifolds, 
in the study of symplectic topology and so on. 

Then in the paper \cite{Ph}, 
a rigorous definition
of the spectral flow was given for not only continuous loops in the space
$\widehat{\mathcal{F}}_*$ but also for arbitrary continuous
paths in $\widehat{\mathcal{F}}_*$ and it was proved that the quantity
is a homotopy invariant of the continuous path in the space
$\widehat{\mathcal{F}}_*$ 
with the fixed end points and satisfies the additivity under catenation
of the continuous paths.  Thus the spectral flow is not only a 
spectral invariant
but also a homotopy invariant, 
in so far as we can treat it in the framework of 
the space $\widehat{\mathcal{F}}_*$.

Nowadays there are many types of
formulas including various ``spectral flows'' corresponding to families of
Fredholm operators which are mostly unbounded operators, because
simply
they
are families of differential operators. 
However in some
cases the ``continuous'' family of such unbounded Fredholm operators 
can be interpreted as a continuous path
in the space $\widehat{\mathcal{F}}_*$ as it was treated in 
the papers \cite{FO}, \cite{BF1} and \cite{CP}. 
Also the theorems in \cite{Fl}, \cite{Yo} and \cite{Ni} 
can be interpreted in the
framework of the space $\widehat{\mathcal{F}}_*$. 
See also \cite{Ta}, \cite{Ge}, \cite{CLM2}, \cite{OF} 
and  \cite{DK} and others.  

It is not clear for the present authors whether arbitrary 
``continuous'' families of Fredholm operators,
especially a family of unbounded selfadjoint
Fredholm operators with varying domains of definitions,
can be interpreted in the framework
of the space $\widehat{\mathcal{F}}_*$.

The main purpose of this paper is to prove a general splitting formula
(Theorem \ref{thm:main}) 
for the spectral flow of a family $\{A+C_t\}_{t\in [0,\,1]}$ 
($C_t$ is of zeroth order) of first order selfadjoint elliptic 
differential operators on a closed manifold in the framework of the
space $\widehat{\mathcal{F}}_*$.
When we decompose the closed manifold $M$ into two parts $M_{\pm}$ 
by a hypersurface $\Sigma$, 
$M = M_{-}\cup M_{+},\, M_{-}\cap M_{+} = \partial M_{\pm} = \Sigma$, 
then we also have a family of
symmetric elliptic operators on each part $M_{\pm}$ simply by the restriction, 
and at this level it is nothing more than having two families of 
symmetric operators defined on manifolds
with boundary.  This means there are no natural choices of boundary
conditions among selfadjoint elliptic boundary conditions 
under which the family 
becomes a family of selfadjoint
operators. Under these circumstances
it will not be apparent whether the spectral flow of the family 
$\{A+C_t\}$ on the
whole manifold is expressed as a sum of two spectral flows obtained by
even imposing a ``suitable'' selfadjoint elliptic boundary 
condition on each part and
adding a correction term 
which solely depends on the boundary manifold $\Sigma$, partly for
the reason that we have no integral representations of the spectral
flow on $M$.

To formulate our splitting formula in this situation,
we shall take the following way, that is, firstly 
we construct a two-parameter family $\{\cal{A}_{s,t}\}_{s,t \in [0,\,1]}$
of selfadjoint Fredholm operators on $M$ with $\cal{A}_{t,t} = A+C_t$.
Then 
by noting the homotopy invariance and the additivity under catenation
of the spectral flow, we will see that 
the sum of two spectral flows of 
$\{\cal{A}_{s,0}\}_{s\in [0,\,1]}$ 
and
$\{\cal{A}_{1,t}\}_{t\in [0,\,1]}$
coincides with the spectral flow of the original family 
(Theorem \ref{thm:pre-main}). 
{}From the definition (see (\ref{equation:two-para})), 
it looks as if the families 
$\{\cal{A}_{s,0}\}_{s\in [0,\,1]}$ 
and
$\{\cal{A}_{1,t}\}_{t\in [0,\,1]}$
are defined on each component respectively;
however in general it is not clear whether each of these two spectral flows
coincides with a spectral flow of a family obtained 
by imposing a suitable selfadjoint elliptic boundary
condition on the restriction of $\{A +C_t\}$ to each part.  
Then secondly
we prove that if the operators in the family are of a product form
near the separating hypersurface, then there is a selfadjoint 
elliptic boundary condition on each part under which the restrictions
of the operators to each part become a family of selfadjoint Fredholm
operators and their spectral flows coincide with each of 
$\{\cal{A}_{s,0}\}_{s\in [0,\,1]}$ 
and 
$\{\cal{A}_{1,t}\}_{t\in  [0,\,1]}$ (Theorem \ref{thm:main}).
Moreover it will turn out that these boundary conditions reflect
the influence from one side to other side of $M_{\pm}$ in a natural way.
We verify these completely in the
framework of the space $\widehat{\mathcal{F}}_*$.  
Our general spectral flow formula (Theorem \ref{thm:general})
and a reduction theorem (Theorem  \ref{thm:reduction})
play a role of a bridge connecting these two spectral flows
for the family of operators of a product form near the separating
hypersurface $\Sigma$. 

There are several similar formulas already treated in the papers
\cite{Ta}, 
\cite{DK}, 
\cite{CLM2} 
and others.  
Here we would like to emphasize that
we can admit non-invertible end points in the family $\{A+C_t\}$ 
because we base our arguments on the rigorous 
definition of spectral flow as a homotopy invariant of the spectral flow given in
\cite{Ph} together with that of the Malsov index which is valid
without any assumptions at the end points
and our method will
explain to some extent what kinds of conditions to the operators in
the family we need to prove such a splitting 
formula. 
Moreover,
since for a family of the operators with product form structure near
the separating hypersurface, the space of our boundary values 
$\boldsymbol\beta = D_{max}/D_{min}$ is identical for any length
of the ``neck'', we expect this property also will give us a reduction 
of our splitting formula (Theorem \ref{thm:general}) under taking
adiabatic limits within our framework.

The contents of this paper are as follows.

In the paper \cite{BF1} a theory was developed 
for the Maslov index in the infinite dimensional case.
However there
we based our arguments on an erroneous characterization
of the Fredholm pair of Lagrangian subspaces 
given in \cite{BW}, which led us to introduce 
a small space of Fredholm-Lagrangian-Grassmannian.
Although 
we corrected this part in the errata 
of the preceding paper \cite{BFO}, here  
in $\S 2$ we will develop the investigation again from the very beginning
of the space of Fredholm-Lagrangian-Grassmannian and the theory of
Maslov index in the infinite dimension. Especially, 
the difference from the well-known finite dimensional case 
will be made clear and we
study a relation of the Maslov index of the continuous curve 
of Fredholm pairs of Lagrangian subspaces with respect to the diagonal
in the product space of symplectic Hilbert spaces and the Maslov index 
of a continuous curve 
with respect to a fixed Lagrangian subspace 
(Proposition \ref{prop:pair}). 

In $\S 3$ we give a new proof of the general spectral flow
formula 
(Theorem \ref{thm:general}) which was proved in the paper \cite{BF1}. 
We prove it
based on the local spectral flow formula (Theorem \ref{thm:local}).
This new proof will give us a clear understanding of the
phenomena surrounding the behavior of the spectrum of selfadjoint Fredholm
operators under a certain perturbation.

In $\S 4$ first we review a reduction theorem (Theorem \ref{thm:reduction})
of the Maslov index in the
infinite dimension which we have proved in the paper
\cite{BFO},
and in the last section $\S 4.3$ 
we prove our main  
Theorem \ref{thm:main}. 
Before proving this Theorem  
we analyze the possibility to express the
spectral flow as a sum of spectral flows of 
two families constructed from the original family
of elliptic operators which are not necessarily of a product
form near the separating hypersurface 
$\Sigma$ (Theorem \ref{thm:pre-main}).  
{}From this formula we arrive at our main Theorem
under the assumption that the operators in the
family are of a product form near the separating hypersurface
$\Sigma$.
The preceding reduction theorem in $\S 4.2$ of the Maslov index is
used twice here, which will explain the role of the product form
assumption on the operators near the separating hypersurface $\Sigma$.

Throughout this paper we shall work in the real category, that is,
Hilbert spaces are real and elliptic operators are defined 
on a real vector bundle.

Finally we remark that although there is a deep theory of 
pseudo-differential operators with transmission property 
including the theory of Calder\'on projector (\cite{Ho2}), 
here we avoid the use of Calder\'on
projector in the $L_2$-framework 
to deal with the Cauchy data space.  This is because
it will make things confusing to use 
the Calder\'on projector from the beginning in the
$L_2$-framework, and furthermore it will force us to miss the role of the
product form assumption on the operators in the splitting formula.
However we must employ two facts without proofs from
the pseudo-differential operator theory:
\begin{enumerate}
\item 
the space of boundary values $\boldsymbol\beta =D_{max}/D_{min}$
of sections
in the maximum domain in $L_2$-space of a
first order selfadjoint elliptic differential operator on a manifold with
boundary $\Sigma$ is included in the
Sobolev space of order $-\frac{1}{2}$ 
on the boundary manifold (\cite{Ho1}).
\item 
if a first order symmetric elliptic differential operator is of 
a product form near the boundary manifold
(= separating hypersurface in our setting), then the closed
extension defined by
the Atiyah-Patodi-Singer boundary condition is selfadjoint and satisfy
``coercive estimate'', which we need to determine the space $\boldsymbol\beta$
(Proposition \ref{prop:4.6}, also see \cite{APS}).
\end{enumerate}


\section{Maslov index in the infinite dimension}

\subsection{Fredholm-Lagrangian-Grassmannian}

Let $(H,\inner{\cdot}{\cdot},\omega)$ be a symplectic 
(real and separable) Hilbert space. Here $\inner{\cdot}{\cdot}$ denotes
the inner product and $\omega$ the symplectic form, i.e., a
non-degenerate, skew-symmetric bounded bi-linear form.

In the discussion below we do not replace the symplectic form $\omega$
after once it was introduced in the real Hilbert space $H$, but we may
always assume that there exists a bounded operator $J:\;H\rightarrow
H$ such that $\omega(x,y)=\inner{Jx}{y}$ for any $x,y\in H$ and
$J^2=-\Id$ by replacing the inner product by another one which defines
an equivalent norm on $H$.  
So we can assume from the beginning the following relations:
\begin{equation*}
\begin{array}{l}
\text{$\trans{J}=-J,\quad \inner{Jx}{Jy}=\inner{x}{y}$\quad and}\\
\omega(Jx,Jy)=\omega(x,y)\quad\text{for all $x,y\in H$.}
\end{array}
\end{equation*}
Here $\trans{J}$ denotes the transpose of $J$ with respect to the
inner product $\inner{\cdot}{\cdot}$.

Let $\cal{L}(H)$ denotes the set of all Lagrangian subspaces of $H$.

Each Lagrangian subspace is closed, and the topology on
$\cal{L}(H)$ is given by embedding it into the space $\cal{B}(H)$ of
bounded linear operators on $H$ by identifying each Lagrangian
subspace with  the orthogonal projection onto it.

We fix a Lagrangian subspace $\lambda\in\cal{L}(H)$.
\begin{defn}\label{def:2.1}
\begin{enumerate}
\item 
The Fredholm Lagrangian Grassmannian of $H$ with respect to $\lambda$
is defined as
\begin{equation*}
\cal{FL}_{\lambda}(H):=\{\mu\in\cal{L}\;|\;(\mu,\lambda)\;
\text{is a Fredholm pair}\}.
\end{equation*}
\item 
We call the subset 
\begin{equation*}
\frak{M}_{\lambda}(H):= \{\mu\in\cal{FL}_{\lambda}(H)\;|\; \mu\cap\lambda\ne
\{0\}\}
\end{equation*}
the Maslov cycle with respect to $\lambda$.
\end{enumerate}
\end{defn}

Note that a pair $(\mu,\lambda)$ of closed subspaces 
is called a Fredholm pair, if it satisfies
\begin{align*}
\text{(1)}\quad &
\dim \mu\cap \lambda < +\infty,\\
\text{(2)}\quad &
\mu +\lambda ~\text{is closed and} ~\dim H/(\mu + \lambda)< +\infty.
\end{align*}
In our case, i.e., for a Fredholm pair $(\mu,\lambda)$ of Lagrangian
subspaces 
$
\dim \lambda \cap\mu = \dim H/(\lambda+\mu)
$.

When we regard the real Hilbert space $H$ as a complex Hilbert
space by the almost complex structure $J$ with the obvious Hermitian
inner product 
$\inner{\cdot}{\cdot}_{J}
=\inner{\cdot}{\cdot}-\sqrt{-1}\omega(\cdot,\cdot)$, we denote it by $H_J$.
Then for a fixed $\lambda \in \cal{L}(H)$ we have
an identification 
\begin{equation}\label{eq:decom}
\begin{array}{l@{}l}
H_{J}=\lambda\oplus\lambda^{\perp}= & \lambda\oplus J\lambda\cong
\lambda\otimes\C\\
& x+Jy\mapsto x\otimes 1+y\otimes\sqrt{-1},
\end{array}
\end{equation}
with the complex conjugate map $\tau_{\lambda}(x+J(y))=x-J(y)$, $x, y
\in \lambda$.

We denote the group of unitary operators on $H_J$ by
\begin{equation*}
\cal{U}(H_J) =\{U\in\cal{B}(H)\;|\; U\,J=J\,U \;\text{and}\; U\,{}^tU = {}^tU\,U=\Id \}.
\end{equation*}

We have a natural action of $\cal{U}(H_J)$ on $ \cal{L}(H)$ and so by
fixing a $\lambda \in  \cal{L}(H)$ we have a surjective map $\rho$:
\begin{equation}\label{eq:n-map}
\begin{array}{r@{\,}l}
\rho:\cal{U}(H_J) & \rightarrow \cal{L}(H)\\
U &\mapsto U(\lambda^{\perp}).
\end{array}
\end{equation}

Put $\cal{U}_{\lambda}(H_J):=\rho^{-1}(\cal{FL}_{\lambda}(H))\subset
\cal{U}(H_J)$. Then each $U\in\cal{U}_{\lambda}(H_{\lambda})$ is
decomposed into $U=X+\sqrt{-1}Y$, under the the identification (\ref{eq:decom}),
such that
\begin{align}
& X,Y:\lambda\rightarrow\lambda\quad\text{bounded},\\
& X\ytrans{Y}=Y\trans{X}\quad\text{and}\quad X\ytrans{X}+Y\trans{Y}=\Id,\\
& X\quad\text{Fredholm operator}
\end{align}

In the same way as for the finite dimensional case(\cite{Le}),
we define
the map $\cal{S}_{\lambda}$
\begin{equation}\label{eq:S-map}
\cal{S}_{\lambda} : \cal{U}_{\lambda}(H_J)\ni U=X+\sqrt{-1}Y\mapsto W_{\mu}:
=(X+\sqrt{-1}Y)(\trans{X}+\sqrt{-1}\ytrans{Y}),
\end{equation}
where $W_{\mu}$ depends only on $\mu=U(\lambda^{\perp})$.
Note that for $U = X+\sqrt{-1}Y$,
 $\trans{X}+\sqrt{-1}\ytrans{Y}$ = 
$\tau_{\lambda}\circ U^{*} \circ \tau_{\lambda}\equiv {\theta}_{\lambda}(U)$.
So we can write $\cal{S}_{\lambda}(\mu) = U\circ\theta_{\lambda}(U)$.

We can easily check that 
$W_{\mu}+\Id
= 2 X(\trans{X}+\sqrt{-1}\ytrans{Y})$ is a Fredholm operator. 
Also we have for 
$x,y\in \lambda$ the point $z=x\otimes 1+y\otimes\sqrt{-1}\in H_J$ 
belongs to $\mu$ if and only if $-z=W_{\mu}(\tau_{\lambda}(z))$, so 
\begin{equation}\label{eq:2.7}
\Ker(W_{\mu}+\Id)=(\mu\cap\lambda)\otimes\C\cong(\mu\cap\lambda)\oplus
J(\mu\cap\lambda).
\end{equation}
Hence
\begin{prop}
For any $\mu\in\cal{FL}_{\lambda}(H)$ and any $U\in\cal{U}_{\lambda}
(H_J)$ with $\mu=U(\lambda^{\perp})$, 
\begin{equation*}
\dim_{\R}(\mu\cap\lambda)=\dim_{\C}\Ker(W_{\mu}+\Id).
\end{equation*}
\end{prop}

Let us now consider the space
\begin{equation}\label{eq:2.8}
\cal{U}_{\cal{F}}(H_J):=\{U\in\cal{U}(H_J)\;|\;\text{$U+\Id$ is a Fredholm
operator}\}
\end{equation}
and a subset
\begin{equation}\label{eq:2.9}
\cal{U}_{\frak{M}}(H_J):=\{U\in\cal{U}_{\cal{F}}(H_J)\;|\;\Ker(U+\Id)\ne
\{0\}\}
\end{equation}
which we can regard as a kind of the universal Maslov cycle, since 
$\cal{S}_{\lambda}^{-1}(\cal{U}_{\frak{M}}(H_J)) =
\frak{M}_{\lambda}(H)$
for any $\lambda$. 

Now we state the fundamental property for discussing the Maslov index in
the infinite dimension.

\begin{prop}\label{prop:pi_1}
\begin{equation*}
\pi_1(\cal{FL}_{\lambda}(H))\simeq\mathbb{Z},\;
\pi_1(\cal{U}_{\cal{F}}(H_J))\simeq\mathbb{Z}
\end{equation*}
and the induced map $(\cal{S}_{\lambda})_{*}:\;\pi_1(\cal{FL}_{\lambda}(H))
\rightarrow\pi_1(\cal{U}_{\cal{F}}(H_J))$ is an isomorphism.
\end{prop}
\begin{proof}
We know already that 
$\pi_1(\cal{FL}_{\lambda}(H))\simeq \mathbb{Z}$.
The fact $\pi_1(\cal{U}_{\cal{F}}(H_J))\simeq\mathbb{Z}$ is proved in
the appendix.
There we will explain a more natural meaning of the space
$\cal{U}_{\cal{F}}(H_J)$ and the proof is similar to that of the first fact
(see also a recent work by Kirk-Lesch \cite{KL}, where a simple proof is
given based on the Calkin algebra).
\end{proof}

\subsection{Maslov index for paths}

In this subsection we review the functional analytic definition of 
the Maslov index.

Let
\begin{equation*}
\begin{array}{rcc}
W:\;I=[0,1]&\rightarrow &\cal{U}_{\cal{F}}(H_J)\\
t&\mapsto & W(t)
\end{array}
\end{equation*}
be a continuous path in $\cal{U}_{\cal{F}}(H_J)$.
Since for Fredholm operators $0$ is an isolated eigenvalue
with finite multiplicity, we have the following property:
\begin{lem}\label{lem:Ph}
There exists a partition $0=t_0<t_1<\dots<t_N=1$ of the interval $I$
and positive numbers 
$\varepsilon_j\,(j=1,\cdots,N)$ with $0<\varepsilon_j<\pi$
such that for each $t$ in $[t_{j-1},\, t_j]$ the operator 
$W(t)-e^{\sqrt{-1}(\pi+\varepsilon_j)}$ is invertible,
and also for each $t$ in  $[t_{j-1},\,t_j]$ the operator $W(t)$ has
only a finite number
of eigenvalues 
in the arc $\{e^{\sqrt{-1}(\pi + \theta)}\,|\, |\,\theta\,|\le\varepsilon_j\}$.
\end{lem}

We now define an `index'
$\mathbf{M}(\{W(t)\}$ 
of the curve $\{W(t)\}_{t\in I}$. 

\begin{defn}\label{def:2.5}
We set
\begin{equation*}
\mathbf{M}(\{W(t)\}):=\sum_{j=1}^N (k(t_j,\varepsilon_j)
-k(t_{j-1},\varepsilon_{j}))
\end{equation*}
with
\begin{equation*}
k(t,\epsilon_j):=\sum_{0\le\theta\le\varepsilon_j}\dim
\Ker(W(t)-e^{\sqrt{-1}(\pi+\theta)})
\end{equation*}
for $t_{j-1}\le t\le t_j$.
\end{defn}

As for Phillips's definition of the spectral flow (see \cite{Ph})
we see that the definition of the `index' does not depend on the
choice of the partition $0 = t_0 \,< t_1\,< \cdots \,< t_N = 1$ of 
the interval and the positive numbers
$\{\varepsilon_j\}_{j=1}^{N}$ satisfying the above Lemma \ref{lem:Ph} .

This `index' has the following properties:
\begin{enumerate}
\item 
additive under catenation of the paths 
\item
Only a homotopy invariant of curves
in $\cal{U}_{\cal{F}}(H_J)$ with fixed endpoints and distinguishes
the homotopy class.
\end{enumerate}
By making use of this `index' we obtain now a functional analytic 
definition of the Maslov index for continuous paths 
in $\cal{FL}_{\lambda}(H)$.  
Let $\mu:\;I\rightarrow\cal{FL}_{\lambda}(H)$ be a continuous path
in $\cal{FL}_{\lambda}(H)$ (so that $\cal{S}_{\lambda}\circ\mu$ is a
continuous path in $\cal{U}_{\cal{F}}(H_J)$).

\begin{defn}\label{def:M-index-path}
We define the Maslov
index of the curve $\{\mu(t)\}$ with respect to $\lambda$ by
\begin{equation*}
\mathbf{Mas}(\{\mu(t)\},\lambda)
:=\mathbf{M}(\{\cal{S}_{\lambda}(\mu(t))\}).
\end{equation*}
\end{defn}

By Proposition \ref{prop:pi_1}, the Maslov index inherits the
properties of the `index'.

\begin{rem}\label{def:almost-coincide}
Let $\lambda '$ be another Lagrangian subspace such that $\lambda\cap
\lambda '$ is of finite codimension in $\lambda$ 
(hence also in $\lambda '$).  We say in this case that $\lambda$ and
$\lambda '$ almost coincide or are almost coincident. This is an
equivalent relation in $\cal{L}(H)$ and we have
$\cal{FL}_{\lambda}(H)$ =
$\cal{FL}_{\lambda '}(H)$ for almost coincident pair $\lambda$ and
$\lambda '$, and the difference of the Maslov index
\begin{equation*}
\mathbf{Mas}(\{\mu(t)\},\lambda) - \mathbf{Mas}(\{\mu(t)\},\lambda ')
\end{equation*}   
only depends on the four Lagrangian subspaces $\mu(0), \mu(1),
\lambda$ and $\lambda '$.  We denoted this quantity by 
$\sigma (\mu(1), \mu(0)\,;\, \lambda, \lambda ')$,
and called it H\"ormander index in the infinite dimension 
as a corresponding quantity to that for the finite dimensional case (\cite{BF2}).
This index works as the transition function
in describing the universal covering space of $\cal{FL}_{\lambda}(H)$.

In fact, let $\widetilde{\cal{FL}_{\lambda}}(H)$ be the universal
covering space of $\cal{FL}_{\lambda}(H)$. We regard this space as the
space of homotopy classes of continuous curves with the common initial point =
$\lambda^{\perp}$ together with its end point. 
{}For $\mu\in \cal{L}(H)$
with the property that $\dim \lambda/(\mu\cap\lambda) < +\infty$,
let us denote $\cal{FL}_{\mu}^{(0)} = \cal{FL}_{\mu}(H)\backslash
\frak{M}_{\mu}(H)$ = $\cal{FL}_{\lambda}(H)\backslash
\frak{M}_{\mu}(H)$
and define a map $g_{\mu}: \cal{FL}_{\mu}^{(0)}\times\Z\rightarrow 
\widetilde{\cal{FL}_{\lambda}}(H)$ by 
\begin{equation*}
g_{\mu}(x,n) = (\text{homotopy class of a curve}\,\, \{\gamma(t)\},\,\gamma(1))
\end{equation*}
where $\{\gamma(t)\}$ is such a curve that $\gamma(0) =
\lambda^{\perp}$, $\gamma(1) = x$ and $\Mas(\{\gamma(t)\},\mu )= n$.
Then the map $g_\mu$ gives a homeomorphism between
$\cal{FL}_{\mu}^{(0)}\times \Z$ and the subspace 
$\pi^{-1}(\cal{FL}^{(0)}_{\mu})$, where 
$\pi :\widetilde{\cal{FL}_{\lambda}}(H)\rightarrow\cal{FL}_{\lambda}(H)$
is the projection map.  Now we know that such subspaces
$\cal{FL}^{(0)}_{\mu}$ cover $\cal{FL}_{\lambda}(H)$ and for $x\in
\cal{FL}^{(0)}_{\mu}\cap
\cal{FL}^{(0)}_{\nu}$, $g_{\nu}^{-1}\circ g_{\mu}(x,n) = (x,
\sigma (x,\lambda^{\perp}; \mu,\nu)+m)$, and we have the additivity $\sigma(x,\lambda^{\perp};\mu,\nu)
+ \sigma(x,\lambda^{\perp}; \nu, \theta) = \sigma(x,\lambda^{\perp};\mu,\theta)$
just by the definition. These show that 
$\cal{FL}^{(0)}_{\mu}\cap\cal{FL}^{(0)}_{\nu}\ni x 
\mapsto \sigma(x,\lambda^{\perp};\mu,\nu)$ is the transition function
for the universal covering space of $\cal{FL}_{\lambda}(H)$.
\end{rem}

\begin{rem}\label{rem:Leray}
In the finite dimensional case,
the Maslov index for continuous paths was first defined in the paper \cite{Go1}
by noting the extendibility of the Leray index for arbitrary pairs of
points on the universal covering space of Lagrangian Grassmannian
through the Kashiwara index and by making
use of the relation between Leray index and Maslov index.  Conversely,
first we define Maslov index for arbitrary paths with
respect to a Maslov cycle as we gave above, then we can define Leray index for
arbitrary pairs of points on the universal covering of Lagrangian Grassmannian.

In the infinite dimensional case we can define the Maslov index for arbitrary
paths with respect to a Maslov cycle as we did above. 
But we cannot define Kashiwara
index for arbitrary triples of Lagrangian subspaces as in the finite
dimensional case, although we have a symmetric bilinear form similar
to the finite dimensional case.  We can define it for mutually almost
coincident triples, since then the symmetric bilinear form is of
finite rank.   Also we cannot define Leray index
for arbitrary pairs of points on the universal covering space of 
Fredholm-Lagrangian-Grassmannian.  It might be possible only for pairs
of unitary operators $U,V$ having the property that $U - V$ is of trace
class through the embedding $\cal{S}_{\lambda}$ (see \cite{KL}).
\end{rem}

For differentiable curves, there is another way of describing
the `index' locally, which is analogous to Robbin-Salamon \cite{RS}.
Let $\{W(t)\}$ be a $C^1$-path in $\cal{U}_{\cal{F}}(H_J)$ with
the differentiability with respect to the operator norm.

\begin{defn}\label{def:2.7}
\begin{enumerate}
\item A parameter $t^*$ with $0<t^*<1$ is called a crossing for the
family $\{W(t)\}$, if $\Ker(W(t^*)+\Id)\ne\{0\}$.
\item We define the crossing form at a crossing $t^*$ as a symmetric
  bilinear form on $\Ker(W(t^*)+\Id)$ by
\begin{equation*}
\tilde{Q}_{\frak{M}}(x,y):=\left.\frac{d}{dt}<x,R(t)y>\right|_{t=t^*}
\quad\text{for $x,y\in\Ker(W(t^*)+\Id)$,}
\end{equation*}
where $\{R(t)\}$ is a family of bounded selfadjoint operators
defined by the condition $W(t)=W(t^*)e^{\sqrt{-1}R(t)}, R(t^*)=0$.
\item A crossing $t^*$ is called regular if the form $\tilde{Q}_{\frak{M}}$
is non-degenerate at $t^*$.
\end{enumerate}
\end{defn}

\begin{prop}\label{prop:M-sign}
Let $\{W(t)\}$ be a path in $\cal{U}_{\cal{F}}(H_J)$ of class $C^1$
and $0<t^*<1$ a regular crossing. Then there exists a real $\delta>0$
such that
\begin{equation*}
\mathbf{M}(\{W(t)\}_{|t-t^*|\le\delta})=\sign\tilde{Q}_{\frak{M}}.
\end{equation*}
\end{prop}
\begin{proof}
Suppose
\begin{equation*}
\dim_{\C}\Ker(W(t^*)-e^{\sqrt{-1}\pi})=k>0,
\end{equation*}
then, since the eigenvalue $-1$ of $W(t^*)$ is isolated, there exists
a positive number $\varepsilon>0$ such that for $0 < |\theta| \le\varepsilon$
the operator $W(t^*)-e^{\sqrt{-1}(\pi+\theta)}$ is invertible. 
Then we can choose a positive number  $\delta_0>0$ such that for $|t-t^*|\le\delta_0$ 
the operator 
$W(t)-e^{\sqrt{-1}(\pi\pm\varepsilon)}$ is invertible, and that
\begin{equation*}
\sum_{|\theta|\le\varepsilon}\dim\Ker(W(t)-e^{\sqrt{-1}(\pi+\theta)})=k,
\end{equation*}
for $|\,t-t^*\,| \,\le\, \delta_0$.
If the signature of $\tilde{Q}_{\frak{M}}$ is $(p,q)$, then by the regularity assumption
at the crossing $t^*$, there exist eigenvalues of
$\dot{R}(t^*)|_{\Ker(W(t^*)+\Id)}$ such that
\begin{equation*}
0<\lambda_1\le\dots\le\lambda_p,\quad 0>\mu_1\ge\dots\ge\mu_q.
\end{equation*}
By Kato's theorem (Theorem VIII 2.6 in \cite{Ka}), there exists a $\delta>0$
(with $\delta<\delta_0$) such that for $t^* \le t\le t^*+\delta$, 
eigenvalues $\{\lambda_i(t)\}_{i=1}^{p}$ and $\{\mu_i(t)\}_{i=1}^q$ 
of $W(t)+\Id$ satisfy the following asymptotic property:
\begin{align*}
\lambda_i(t)&=\exp(\sqrt{-1}(\pi+\lambda_i(t-t^*)+o(|t-t^*|)))\quad(i=1,\cdots,p),\\
\mu_j(t)&=\exp(\sqrt{-1}(\pi+\mu_j(t-t^*)+o(|t-t^*|)))\quad(j=1,\cdots,q).
\end{align*}
Therefore we have
\begin{align*}
&\sum_{0\le\theta\le\varepsilon}\dim\Ker(W(t)-e^{\sqrt{-1}(\pi+\theta)})=p,\\
&\sum_{-\varepsilon\le\theta<0}\dim\Ker(W(t)-e^{\sqrt{-1}(\pi+\theta)})=q
\end{align*}
for $t^* < t\le t^*+\delta$. For $t^*-\delta\le t<t^*$, $p$ and
$q$ are exchanged.

Hence
\begin{align*}
\mathbf{M}(\{W(t)\}_{|t-t^*|\le\delta})
&= k(t^*+\delta,\varepsilon)-k(t^*-\delta,\varepsilon)\\
&=p-q=\sign\tilde{Q}_{\frak{M}}.
\end{align*}
\end{proof}

\begin{rem}
\begin{enumerate}
\item Kato's theorem (\cite{Ka}) is proved 
for analytic curves, but this is also
true for $C^1$-curves.
\item For crossing $t^*=0$ or $1$, we only consider the
one-sided differentiation in the definition of the crossing form.
In these cases we have
\begin{align*}
\mathbf{M}(\{W(t)\}_{0\le t\le \delta})&=-q,\\
\mathbf{M}(\{W(t)\}_{1-\delta\le t\le 1})&=p',
\end{align*}
where the signature of $\tilde{Q}_{\frak{M}}$ at $t^*=0$ is $(p,q)$ and
at $t^*=1$\, $(p',q')$.
\end{enumerate}
\end{rem}

\begin{cor}\label{cor:2.11}
Let $\mu:\;I\rightarrow\cal{FL}_{\lambda}(H)$ be a $C^1$-class path
(so that $\cal{S}_{\lambda}\circ\mu(t)$ is a path in
$\cal{U}_{\cal{F}}(H_J)$ also of class $C^1$). Let
$0<t^*<1$ be a regular crossing of the curve 
$\{\cal{S}_{\lambda}\circ\mu(t)\}$.
Then there exists a $\delta>0$ such that
\begin{equation*}
\mathbf{Mas}(\{\mu(t)\}_{|t-t^*|\le\delta},\lambda)=\sign\tilde{Q}_{\frak{M}},
\end{equation*}
where $\tilde{Q}_{\frak{M}}$ denotes the crossing form of $\{\cal{S}_{\lambda}\circ
\mu(t)\}$ at the time $t=t^*$.
\end{cor}

There is another differential description of the Maslov index
which will turn out to be more suitable for proving the spectral
flow formula. It is based on a representation of $\mu$ as the graph
of a suitable bounded operator. Let $\mu:\;I\rightarrow
\cal{FL}_{\lambda}(H)$ be a path in $\cal{FL}_{\lambda}(H)$ of class
$C^1$ and let $0<t^*<1$ be a crossing of the curve $\{\cal{S}_{\lambda}\circ
\mu(t)\}$, i.e., $\mu(t^*)\cap\lambda\ne\{0\}$. Put $\mu=\mu(t^*)$,
then $\mu(t)$ is transversal to $\mu^{\perp}$ for $|t-t^*|\ll 1$
and in this neighborhood of $t^*$, each $\mu(t)$ can be written
as the graph of the bounded operator $\varphi(t):\;
\mu\rightarrow\mu^{\perp}$. Note that the curve $\{\varphi(t)\}$
is also of class $C^1$. We consider the bilinear form
\begin{equation*}
Q_{\frak{M}}(x,y):=\left.\frac{d}{dt}\omega(x,\varphi(t)y)\right|_{t=t^*}
\quad\text{for $x,y\in\mu(t^*)\cap\lambda$.}
\end{equation*}
It is a symmetric bilinear form and we have the following proposition.

\begin{prop}
If $\mu(t^*)\cap\lambda\ne\{0\}$, then $\sign Q_{\frak{M}}=
\sign\tilde{Q}_{\frak{M}}$.
\end{prop}

\begin{proof} 
  For $|t-t^*|\ll 1$, $\mu(t)$ is represented in two ways;
\[  \mu(t)=U(t)(\lambda^\bot)=\mbox{graph of}\, \varphi(t),
\]
where $U(t)$ is a unitary operator $\in \cal{U}_{\lambda}(H_J)$ and 
$\varphi :\mu\to\mu^\bot,\quad \mu=\mu(t^*),\quad \varphi(t^*)=0$.

  We write operators $U(t)$ and $W(t)=U(t)\circ \theta_{\lambda}(U(t))$
in the following form;
\[
   U(t)=U(t^*)e^{\sqrt{-1}A(t)}, \quad W(t)=W(t^*)e^{\sqrt{-1}R(t)},
\]
where $A(t)$ and $R(t)$ are selfadjoint operators
with $A(t^*)=0$ and $R(t^*)=0$.

Now put $A(t)=X(t)+\sqrt{-1}Y(t)$ with $X= {^tX}$, $Y=-{^tY}$,
then we obtain relations among differentials
of $A(t)$, $R(t)$ and $\varphi(t)$ at $t=t^*$ as follows;
\begin{align*}
\theta_{\lambda}(U(t^*))\dot{R}(t^*) &
= 2\dot{X}(t^*)\circ\theta_{\lambda}(U(t^*)) \\
U(t^*)^{-1}(-J\dot{\varphi}(t^*)) &
= \dot{X}(t^*)U(t^*)^{-1}
\quad\mbox{on}\,\mu\cap\lambda .
\end{align*}   

 Therefore $\dot{R}(t^*)$ and $-J\dot{\varphi}(t^*)$
have the same signature on $\mu\cap\lambda$.
By recalling that signatures of $R(t^*)$ and $-J\dot{\varphi}(t^*)$
are $\sign {\widetilde Q}_{\frak{M}}$ and $\sign Q_{\frak{M}}$ respectively,
 we complete the proof.
\end{proof}


\subsection{Maslov index for the curve of Fredholm pairs of
Lagrangian subspaces}

In the following we will denote the direct sum of the symplectic
Hilbert space $(H,\omega)$ and $(H, -\omega)$ with the notation 
$\Hi=H\boxplus H$ $\equiv$ $H_{\omega}\oplus H_{-\omega}$. 
$\Hi$ is a symplectic Hilbert space with the symplectic form
$\Omega(x\boxplus y,~u\boxplus v) = \omega(x,u)~-~\omega(y,v)$, 
and the corresponding almost complex structure
$\J = J\oplus -J$, so that we have 
$\Hi_{\J}=H_J\oplus H_{-J}$.
 
Let $\{(\mu_t,\lambda_t)\}_{t\in I}$ be a continuous family of
Fredholm pairs of Lagrangian subspaces, then 
$\{\mu_t\boxplus\lambda_t\}$
is a curve in $\cal{FL}_{\Delta}(H\boxplus H)$,
where $\Delta$ is the diagonal of $H\boxplus H$. 
Of course it is natural to define the Maslov index of the curve
of Fredholm pairs $\{(\mu_t,\lambda_t)\}$ to be
$\Mas(\{\mu_t\boxplus\lambda_t\},\Delta)$. 

\begin{prop}\label{prop:pair}
Let $\{\mu_t\}$ be a continuous curve in $\cal{FL}_{\lambda}
(H_{\omega})$, then
\begin{equation*}
\Mas(\{\mu_t\},\lambda)=\Mas(\{\mu_t\boxplus\lambda\},\Delta).
\end{equation*}
\end{prop}

\begin{rem}
For loops the property  will be well-known. For arbitrary
continuous paths, in the finite dimensional case 
this may be proved by making use of
relations with Leray index as we noted in Remark \ref{rem:Leray}, but 
in the infinite dimensional case we have no such relations and we 
need a proof which is valid not only for loops but also
for any continuous paths.
\end{rem}

If we identify 
$\Hi_{\J}= \Delta+\Delta^{\perp}
=\Delta+\J(\Delta)\cong\Delta\otimes \C$, then 
$\tau_{\Delta}(a\boxplus b)=b\boxplus a$.
Let us decompose $H$ as $H= \lambda\oplus\lambda^{\perp}$ and let $\varphi:\;\Delta
\rightarrow\Delta^{\perp}$  be
\begin{equation*}
\varphi((x,y)\boxplus(x,y))=(-x,y)\boxplus(x,-y)
\end{equation*}
where we express elements in $\Delta$ by $(x,y)\boxplus (x,y)$,
\,$x+y\in\lambda+\lambda^{\perp}=H$. Then we have
\begin{equation*}
\graph\varphi=\lambda^{\perp}\boxplus\lambda.
\end{equation*}

Let $A=\J\circ\varphi:\;\Delta\rightarrow\Delta$ and $V:\;\Hi_{\J}
\rightarrow\Hi_{\J}$ by
\begin{equation*}
V=\frac{-\sqrt{-1}}{\sqrt{2}}-\frac{A\otimes\Id}{\sqrt{2}}
\end{equation*}
where we regard $A=A\otimes Id$ is complexified according 
to the identification $\Hi_J\cong\Delta\otimes\C$. Then we have
\begin{equation}
\sqrt{-1}(A\otimes\Id)((a,b)\boxplus (c,d)) = (c,-d)\boxplus(-a,b)
\end{equation}
for $(a,b)\boxplus (c,d)\in H_J\boxplus H_{J} =
(\lambda+\lambda^{\perp})\boxplus (\lambda+\lambda^{\perp})$ 
and
\begin{equation}
V(\Delta^{\perp})=\lambda^{\perp}\boxplus\lambda.
\end{equation}

Now we define maps ${\bf{a}}_{\lambda}, \,{\bf{b}}_{\lambda}$ and
$P_{\lambda}$ as follows:
\begin{equation*}
\begin{array}{rcc}
{\bf{a}}_{\lambda}:\;\cal{U}_{\lambda}(H_J)&\longrightarrow&\cal{U}_{\Delta}(\Hi_{\J})\\
U & \mapsto & \tilde{U}\circ V
\end{array}
\end{equation*}
where $\tilde{U}=U\oplus\Id\,:\,H_J\oplus H_{-J} \longrightarrow
H_J\oplus H_{-J}$,
\begin{equation*}
\begin{array}{rcc}
{\bf{b}}_{\lambda}:\;\cal{U}_{\cal{F}}(H_J)&\longrightarrow &\cal{U}_{\cal{F}}
(\Hi_{\J})\\
W & \mapsto & \sqrt{-1}\cdot W\circ (A\otimes\Id),
\end{array}
\end{equation*}
and
\begin{equation*}
\begin{array}{rcc}
P_{\lambda}\,:\,\cal{FL}_{\lambda}(H)&\rightarrow&\cal{FL}_{\Delta}(\Hi)\\
\mu & \mapsto & \mu\boxplus\lambda.
\end{array}
\end{equation*}

\begin{lem}
The following diagram is commutative.
\begin{equation*}
\begin{CD}
\cal{U}_{\lambda}(H_J) @>{{\bf{a}}_{\lambda}}>>
\cal{U}_{\Delta}(H_J\oplus H_{-J})\\
@V{\rho_{\lambda}}VV  @VV{\rho_{\Delta}}V\\
\cal{FL}_{\lambda}(H) @>>{P_{\lambda}}>
\cal{FL}_{\Delta}(H\boxplus H)\\
@V{\cal{S}_{\lambda}}VV @VV{\cal{S}_{\Delta}}V\\
\cal{U}_{\cal{F}}(H_J)@>{{\bf{b}}_{\lambda}}>>\cal{U}_{\cal{F}}(H_J
\oplus H_{-J})
\end{CD}
\end{equation*}
\end{lem}

\begin{proof}
It will be enough to prove 
$\cal{S}_{\Delta}\circ P_{\lambda}={\bf{b}}_{\lambda}\circ
\cal{S}_{\lambda}$. Since  
$\theta_{\Delta}(V) = V$, $V^2 = \sqrt{-1}\cdot A\otimes\Id$ and
$\theta_{\Delta}({\widetilde{U}}) = \Id\boxplus U^{*}$ we have
\begin{eqnarray*}
&&\cal{S}_{\Delta}\circ\rho_{\Delta}({\bf{a}}_{\lambda}(U))\\
&=&\tilde{U}\circ V\circ\tau_{\Delta}\circ(\tilde{U}\circ V)^*\circ
\tau_{\Delta}\\
&=&\tilde{U}\circ\sqrt{-1}(A\otimes\Id)\circ\theta_{\Delta}(\tilde{U})\\
&=&{\widetilde{U\circ\theta_{\lambda}(U)}}\circ \sqrt{-1}(A\otimes\Id),
\end{eqnarray*}
which prove the commutativity of the diagram.
\end{proof}

\renewcommand{\proofname}{Proof of Proposition \ref{prop:pair}}
\begin{proof}
{}From the above lemma we can show that if $E$ is an eigenvalue of
$\cal{S}_{\Delta}(\rho_{\Delta}({\bf{a}}_{\lambda}(U)))$, then $-E^2$
is an eigenvalue of $\cal{S}_{\lambda}\circ\rho_{\lambda}(U)$.
Conversely if $l=e^{\sqrt{-1}\sigma}$ is an eigenvalue of $\cal{S}_{\lambda}\circ
\rho_{\lambda}(U)$, then only one of $\pm e^{\sqrt{-1}(\pi+\sigma)}$
is close to $-1$. So if we have 
a continuous curve $\{\mu_t\} \subset \cal{FL}_{\lambda}(H)$,
then the numbers of eigenvalues of $\{\cal{S}_{\lambda}(\mu_t)\}$ and
$\{\cal{S}_{\Delta}(\mu_t\boxplus\lambda)\}$ across
$e^{\sqrt{-1}\pi}$ 
coincide in both directions. 
This proves the proposition.
\end{proof}

Also we can prove the following property in a similar way as above:

\begin{prop}
$\Mas(\{\mu_t\boxplus\lambda_t\},\Delta)=-\Mas(\{\lambda_t\boxplus
\mu_t\},\Delta)$.
\end{prop}

\renewcommand{\proofname}{Proof}


\section{Spectral flow formula}

\subsection{Symmetric operators and boundary value space}

Let $H$ be a real separable
Hilbert space and $A$ 
a densely defined closed symmetric operator
with the domain $D_m$. 
Let $A^*$ denote its adjoint with the domain $D_M$. 

Let $\boldsymbol\beta$ be the factor space of $D_M$ by $D_m$, $\boldsymbol\beta$ =
$D_M/D_m$,
and let $\gamma$ $:$ $D_M\rightarrow \boldsymbol\beta$ $(\gamma(x) = [x])$ 
be the projection map.
The space $\boldsymbol\beta$ becomes a symplectic Hilbert space with the inner
product induced by the graph norm 
\begin{equation*}
\inner{x}{y}_G:=\inner{x}{y}+\inner{A^*x}{A^*y}
\end{equation*}
on $D_M$
and the symplectic form given by `Green's form'
\begin{equation*}
\omega([x],[y]):=\inner{A^*x}{y}-\inner{x}{A^*y}\quad
\text{for $x,y\in \boldsymbol\beta$.}
\end{equation*}
For $D_m\subset D\subset D_M$, we denote 
\begin{equation*}
A_D:=A^*|_D
\end{equation*}
and the following relations are fundamental for our theory: 
\begin{center}
\begin{tabular}{rcl}
$A_D$ is closed & $\Leftrightarrow$ & $\gamma(D)$ is closed\\
$A_D$ is self-adjoint & $\Leftrightarrow$ & $\gamma(D)$ is Lagrangian\\
$A_D$ has compact resolvent & $\Leftrightarrow$ & 
$D\hookrightarrow H$ is compact in graph norm.
\end{tabular}
\end{center}
(See \cite{DS}, Chapter XII.)

\begin{defn}
We call $\gamma(\Ker A^*)$ the Cauchy data space.
\end{defn}
In fact if $A$ is an elliptic differential operator, then $\boldsymbol\beta$
can be embedded to the distribution space on the boundary
manifold and $\gamma(\Ker A^*)$ is realized as the boundary values
of $L_2$-solutions.

Now we make two assumptions.
\begin{ass}\label{ass:A}
We assume that $A$ admits at least one self-adjoint Fredholm extension
$A_D$ with the domain of definition $D \subset D_M$ 
and $A_D$ has the compact resolvent.
\end{ass}

Then it follows (\cite{BF1}, Proposition 3.5) 
that $\gamma(\Ker (A^*))$ is a
Lagrangian subspace of 
$\boldsymbol\beta$ and $(\gamma(\Ker A^*),\gamma(D))$ is a
Fredholm pair. 

\begin{ass}\label{ass:B}
We assume that there exists a continuous curve $\{C_t\}$ in the
space of bounded self-adjoint operators on $H$ and that the
operators $A^*+C_t-s$ for small $s$ 
satisfy the ``abstract unique continuation property'';
\begin{equation*}
\Ker(A^*+C_t-s)\cap D_m=\{0\}\quad\text{for $|s|\ll 1$ 
in a neighborhood
of each $t$.}
\end{equation*}
\end{ass}

It follows (see \cite{BF1}, Theorem 3.8) that 
$\{\gamma(\Ker(A^*+C_t))\}$
is a continuous curve in $\cal{FL}_{\gamma(D)}(\boldsymbol\beta)$.

\subsection{Definition of the spectral flow}

Let $H$ be a real separable Hilbert space and let $\hat{\cal{F}}$
denote the space of bounded self-adjoint Fredholm operators on $H$.
Following \cite {Ph} 
we recall the definition of 
the spectral flow for continuous paths in $\hat{\cal{F}}$.

Let $\{A_t\}_{t\in[0,1]}$ be a continuous path in $\hat{\cal{F}}$.

\begin{defn}
We define the spectral flow by
\begin{equation*}
\Sf(\{A_t\}):=\sum_{j=1}^N 
(k(t_j,\varepsilon_j)-k(t_{j-1},\varepsilon_j)),
\end{equation*}
where $\{t_j\}_{j=0}^{N}$, a partition of the interval $I$, 
and positive numbers $\varepsilon_j>0$ ($j=1,\cdots,N$) 
have the following properties:
the operator 
$A_t-\varepsilon_j$ is invertible for $t\in [t_{j-1},\, t_j]$ 
and the operator $A_t$ has only a finite number of eigenvalues
in the interval $[-\varepsilon_j,\,\varepsilon_j]$ for each
$t\in [t_{j-1},\,t_j]$.
Here we denote
\begin{equation*}
k(t,\varepsilon_j):=\sum_{0\le\theta<\varepsilon_j}\dim\Ker(A_t-\theta)
\quad\text{for $t_{j-1}\le t\le t_j$.}
\end{equation*}
\end{defn}

Similarly to the case of the Maslov index, the spectral flow is well defined
for homotopy class of curves with fixed endpoints and it
distinguishes the homotopy class of curves in $\hat{\cal{F}}$ and is additive
under catenation of curves. Especially, if the curve is in the trivial
component  $\hat{\cal{F}}_{\pm}$ of  $\hat{\cal{F}}$, 
then it depends only on the end points and equals the difference
$\sum_{\lambda <0} \dim \Ker (A_0-\lambda) - \sum_{\lambda < 0}\dim \Ker
(A_1-\lambda)$,
in case $\{A_t\} \subset \hat{\cal{F}}_+$ and we  have a 
similar formula for $\{A_t\}\subset \hat{\cal{F}}_-$.

For $C^1$-class curves, there is another description of the spectral
flow locally. On $\Ker A_{t^*}\,(0<t^*<1)$, we consider the
well-defined symmetric bilinear form (crossing form) which is
similar to \S{}2.2:
\begin{equation*}
Q_{\Sf}(x,y):=\left.\frac{d}{dt}\inner{x}{A_t(y)}
\right|_{t=t^*}\quad\text{for $x,y\in\Ker A_{t^*}$.}
\end{equation*}
We have an analogous proposition to \ref{prop:M-sign}.

\begin{prop}\label{prop:3.3}
If the quadratic form $Q_{\Sf}$ is non-degenerate on $\Ker A_{t^*}$,
then for sufficiently small $\delta>0$, we have
\begin{equation*}
\Sf(\{A_t\}_{|t-t^*|\le\delta})=\sign (Q_{\Sf}|_{\Ker A_{t^*}}).
\end{equation*}
\end{prop}
\begin{proof}
This is proved in a quite similar way as 
the proof of Proposition \ref{prop:M-sign}.
\end{proof}

Next we proceed to families of unbounded Fredholm self-adjoint
operators $A_t=A_D+C_t$, as discussed in the preceding paragraph.

We assume the continuity of the family $\{C_t\}$ in $\cal{B}(H)$.
First we apply the transformation $A\mapsto\cal{R}(A):=A\sqrt{\Id+A^2}^{-1}$
from the space $C\hat{\cal{F}}$ of selfadjoint Fredholm operators to 
$\hat{\cal{F}}$.
The continuity of maps $t \mapsto \cal{R}(A_D+C_t)$ was proved 
in \cite{BF1} (see also \cite{CP}).

\begin{defn}
We define the spectral flow for $\{A_D+C_t\}_{t\in I}$ as
\begin{equation*}
\Sf(\{A_D+C_t\}):=\Sf(\{\cal{R}(A_D+C_t)\}).
\end{equation*}
\end{defn}

If the family $\{C_t\}$ is $C^1$-class in $\cal{B}(H)$, then we can
define the crossing form as follows;
\begin{equation*}
Q_0(x,y):=\left.\frac{d}{dt}\inner{x}{C_t(y)}\right|_{t=t^{*}}
\quad\text{for $x,y\in\Ker(A^*+C_{t^*})$.}
\end{equation*}

By the proposition 2.10 in \cite{CP}, we know that $\cal{R}(A_D+C_t)$
is a $C^1$-class curve in $\hat{\cal{F}}$ and the derivative is given by
the formula
\begin{align*}
&\frac{d}{dt}(\cal{R}(A_D+C_t))\\
&= \frac{1}{\pi}\int_{0}^{\infty}
s^{-1/2}\{(1+(A_D+C_t)^2+s)^{-1}(1+s)\frac{d}{dt}(C_t)(1+(A_D+C_t)^2+s)^{-1}\\
&-(A_D+C_t)(1+(A_D+C_t)^2+s)^{-1}\frac{d}{dt}(C_t)
(A_D+C_t)(1+(A_D+C_t)^2+s)^{-1}\}ds.
\end{align*}

{}In particular for $x\in \Ker(A_D+C_t)$ we have 
$\frac{d}{dt}(\cal{R}(A_D+C_t))(x)$ = $\frac{d}{dt}C_t(x)$, so we have
\begin{prop}
$Q_0(x,y)=Q_{\Sf}(x,y)$ at crossing $t=t^*$ for any $x, y \in \Ker (A_D+C_{t^*})$.
\end{prop}

\subsection{Local spectral flow formula}

It would be interesting for us to derive a local identification of the spectral
flow with the Maslov index for curves of $C^1$-class. We inherit
the notations and the assumptions from \S{}3.1. We have first the
following
\begin{prop}\label{prop:3.6}
{}For any $t^*\in [0,1]$ and $x,y\in \Ker (A^*+C_t)$,
\begin{equation*}
Q_0(x,y)=Q_{\frak{M}}(\gamma(x),\gamma(y)).
\end{equation*}
\end{prop}
\begin{proof}
When $S_{t^*}$ = $\Ker (A^*+C_t)\ne 0$, write $\gamma(\Ker (A^*+C_t))=
\lambda_t=\graph\varphi_t$, where
$\varphi :\;\lambda_{t^*}\rightarrow\lambda_{t^*}^{\perp}$ for
$|t-t^*|\ll 1$. We can obtain the lift of $\varphi_t$ as follows:
Fix a closed direct summand 
$\mu$ in $D_M$ as $D_M=S_{t^*}\oplus D_m\oplus\mu$, then
we can define
\begin{equation*}
\begin{array}{rcl}
f_t\oplus\psi_t: & S_{t^*} \rightarrow & D_m\oplus\mu\\
x & \mapsto & f_t(x)\oplus\psi_t(x)
\end{array}
\end{equation*}
such that $S_t=\graph f_t\oplus\psi_t$. By the definition of the
crossing form and the symplectic structure, we have
\begin{align*}
Q_{\frak{M}}(\gamma(x),\gamma(y))
&=Q_{\frak{M}}([x],[y])=\left.\frac{d}{dt}\omega([x],\varphi_t[y])
\right|_{t=t^*}\\
&=\left.\frac{d}{dt}\{\inner{A^*x}{\psi_t(y)}-\inner{x}{A^*(\psi_t(y))}
\}\right|_{t=t^*}
\tag{$*$}\label{eq:P3.8-1}
\end{align*}
Since
\begin{equation*}
(A^*+C_t)(y+f_t(y)+\psi_t(y))=0\quad\text{for all $y\in S_{t^*}$,}
\end{equation*}
by differentiating the above equation with respect to $t$ and
evaluating at $t=t^*$, we obtain
\begin{equation*}
\dot{C}_{t^*}(y)+(A^*+C_{t^*})(\dot{f}_{t^*}(y)+\dot{\psi}_{t^*}(y))=0.
\end{equation*}
Therefore,
\begin{equation*}
A^*(\dot{\psi}_{t^*}(y))=-\dot{C}_{t^*}(y)-(A^*+C_{t^*})(\dot{f}_{t^*}(y))
-C_{t^*}(\dot{\psi}_{t^*}(y)).
\end{equation*}
By substituting this result into (\ref{eq:P3.8-1}), we get
\begin{equation*}
\text{(\ref{eq:P3.8-1})}=\inner{A^*x}{\dot{\psi}_{t^*}(y)}
+\inner{x}{\dot{C}_{t^*}(y)+(A^*+C_{t^*})(\dot{f}_{t^*}(y))
+C_{t^*}(\dot{\psi}_{t^*}(y))}.
\end{equation*}
Since $(A^*+C_{t^*})x=0$ and $\dot{f}_{t^*}(y)\in D_m$, we have
\begin{eqnarray*}
&&\inner{A^*x}{\dot{\psi}_{t^*}(y)}+\inner{x}{C_{t^*}(\dot{\psi}_{t^*}(y))}
=\inner{(A^*+C_{t^*})(x)}{\dot{\psi}_{t^*}(y)}=0\\
&&\inner{x}{(A^*+C_{t^*})(\dot{f}_{t^*}(y))}=\inner{x}{(A+C_{t^*})
(\dot{f}_{t^*}(y))}\\
&=&\inner{(A^*+C_{t^*})(x)}{\dot{f}_{t^*}(y)}=0.
\end{eqnarray*}
Therefore we obtain $\text{(\ref{eq:P3.8-1})}=\inner{x}{\dot{C}_{t^*}(y)}$,
which completes the proof.
\end{proof}

Combining Proposition \ref{prop:3.3} with Proposition \ref{prop:3.6}
and Corollary \ref{cor:2.11} yields the following local spectral
formula.

\begin{thm}\label{thm:local}
If the quadratic form $Q_{\Sf}$ is non-degenerate on $S_{t^*}\cap D$,
then we have
\begin{equation*}
\Sf(\{A_D+C_t\}_{|t-t^*|\le\delta})=
\Mas(\{\lambda_t\}_{|t-t^*|\le\delta},\gamma(D))
\end{equation*}
for sufficient small $\delta > 0$.
\end{thm}

\subsection{General spectral flow formula}

Under assumptions made in $\S\, 3.1$ (Assumption (A) and Assumption (B)),
we obtain a general spectral flow formula. This theorem was proved
in \cite{BF1}, Theorem 5.1. 
Here we give a proof of the formula by emphasizing
the local spectral flow formula (Theorem \ref{thm:local}) 
in the preceding section.

\begin{thm}\label{thm:general}
Let $\{A_D+C_t\}$ be a family satisfying Assumption (A) and (B).
Then we have
\begin{equation*}
\Sf(\{A_D+C_t\})=\Mas(\{\gamma(S_t)\},\gamma(D)).
\end{equation*}
\end{thm}

\begin{proof}
Recalling the additivity of the spectral flow and the Maslov index
under catenation of two curves, it is sufficient to prove that
\begin{multline*}
\Sf(\{A_D+C_t\}_{t_1\le t\le t_2})=
\Mas(\{\gamma(S_t)\}_{t_1\le t\le t_2}\},\gamma(D))\\
\text{for $0<t_2-t_1\ll 1$.}
\end{multline*}
We deform the curve $\{A_D+C_t\}_{t_1\le t\le t_2}$ to
$\{A_D+C_{t_1}+s\}_{0\le s\le\varepsilon}\cup\{A_D+C_t+
\varepsilon\}_{t_1\le t\le t_2}\cup\{A_D+C_{t_2}+
\varepsilon-s\}_{0\le s\le\varepsilon}$, provided that for
$t_1\le t\le t_2$ and $0\le s\le\varepsilon$, Assumption (B) is
satisfied.

By the property of homotopy invariance, the spectral flow and the Maslov
index are equal respectively for previous two curves.

Since $\Ker(A_D+C_t+\varepsilon)=0$ for $t_1\le t\le t_2$, the
spectral flow and the Maslov index vanish for corresponding
curve. By applying local spectral flow formula (Theorem \ref{thm:local})
to the smooth curves $\{A_D+C_{t_1}+s\}_{0\le s\le\varepsilon}$
and $\{A_D+C_{t_2}+\varepsilon-s\}_{0\le s\le\varepsilon}$, we get
the result of this theorem.
\end{proof}


\section{Splitting formula for a spectral flow}

In this section we prove our main Theorem \ref{thm:main}.

Let $\mathbb{E}$ be a real vector bundle  
on a closed manifold $M$ (we assume $M$ is connected), and 
let $A$ be an elliptic
and selfadjoint first order differential operator acting on the vector
bundle $\mathbb{E}$. The selfadjointness of the operator $A$ will be
with respect to a suitably fixed $L_2$-inner product on the space of smooth sections 
$C^{\infty}(M,\mathbb{E})$. 
We will not describe the inner product on $\mathbb{E}$ explicitly.
Let $\{C_t\}$ be a continuous family of selfadjoint bundle maps
of $\mathbb{E}$, and we will regard them, using the same notation, as a continuous
family of zeroth order operators on $L_2(M,\mathbb{E})$.  We also denote
by $H^s(M,\mathbb{E})$ a Sobolev space of sections of $\mathbb{E}$ of order $s$.
Throughout this section we assume that the operators 
$A+C_t + s$ ($|s|\ll 1$) 
satisfy the unique continuation property 
for any hypersurface (Assumption (B)).

\subsection{Splitting of manifolds and first order
differential operators}

Let $\Sigma$ be a hypersurface on $M$ and let $D_{min}$ be a subspace
in $H^1(M,\mathbb{E})$ consisting of sections which vanish on $\Sigma$.
Then the operator $A$ defined on $D_{min}$, denoted by
$\cal{A}_0$, is a closed symmetric operator.  We will denote the
domain of the adjoint operator $(\cal{A}_0)^*$ by $D_{max}$.
Here in addition if we assume that the hypersurface separates the manifold $M$
into two parts $M_{\pm}$ with common boundary 
$\Sigma = \partial{M_{\pm}}$, then $D_{min}$ and $D_{max}$ are also decomposed
into two components corresponding to submanifolds $M_{\pm}$.
We denote them by $D_{min} = D_{min}^{-}\oplus D_{min}^{+}$,
$D_{max} = D_{max}^{-}\oplus D_{max}^{+}$ and 
$\cal{A}_0 = \cal{A}_0^{-}\oplus\cal{A}_0^{+}$,
and $D_{max}^{\pm}$ is the domain of the adjoint operator
of $\cal{A}_0^{\pm}$ respectively.  
Although graph inner products on $D_{max}$ by operator
$A+C_t$ are not identical, the norms are all equivalent and
the symplectic forms $\omega_{\pm}$ 
on $\boldsymbol\beta^{\pm} = D_{max}^{\pm}/D_{min}^{\pm}$ 
do not depend on the parameter.  
 
Notice that if $\Sigma$ is only orientable and
does not separate $M$, that is, if $M\backslash\Sigma$ is connected,
then $D_{min}($and also $D_{max}$) does not
decompose in the above way, but $\boldsymbol\beta = D_{max}/D_{min}$
is a sum of two spaces $\boldsymbol\beta^{\pm}$ of boundary values taken
from each side.  However in this case the Cauchy data space 
$\gamma(\Ker(\cal{A}_0)^*)$ is not a sum of two Lagrangian subspaces 
in $\mathbb\beta^{\pm}$.

Henceforth we assume that a hypersurface $\Sigma$ separates $M$
into two parts $M_{\pm}$.
 
\begin{prop}[\cite{Ho1}]
$\boldsymbol\beta^{\pm}$ is a subspace in $H^{-1/2}(\Sigma,\mathbb{E}|_{\Sigma})$, the
Sobolev space of sections with values in the bundle $\mathbb{E}$ restricted to
$\Sigma$, of order $-1/2$. 
Also we have $\boldsymbol\beta^{-}\cap\boldsymbol\beta^{+}
=H^{1/2}(\Sigma,\mathbb{E}|_{\Sigma})$.
\end{prop}

Let $\gamma_{\pm}$ be the map from $D_{max}^{\pm}$ to $\boldsymbol\beta^{\pm}$.
The families
$\{\Lambda_t^{\pm}\}=\{\gamma_{\pm}(\Ker((\cal{A}_0^{\pm}+C_t)^*))\}$  
are continuous families of Lagrangian subspaces of $\boldsymbol\beta^{\pm}$ and
according to the general spectral flow
formula in $\S\, 3$ we have

\begin{thm}  From the ellipticity  of the operator we have 
$\Lambda_t^{-}\oplus\Lambda_t^{+}\in\cal{FL}_{\boldsymbol\delta}(\boldsymbol\beta)$,
and 
$\Sf\{A+C_t\}=\Mas(\{\Lambda_t^{-}\oplus\Lambda_t^{+}\},\boldsymbol\delta)$
where
\begin{equation*}
\boldsymbol\delta =\{(\varphi,\varphi)\in \boldsymbol\beta^{-}\oplus \boldsymbol\beta^{+}\;|\;{}^{\exists}
f\in H^1(M,\mathbb{E}),\gamma_{\pm}(f|_{M_{\pm}})=\varphi\}.
\end{equation*}
\end{thm}

Let $\{L_{s,t}\}$ be a family of bounded operators on $L_2(M,\mathbb{E})$
defined as
\begin{equation*}
L_{s,t}(f)=\left\{\begin{array}{ll}
0 & \text{on $M_{-}$}\\
(C_t-C_s)(f) & \text{on $M_{+}$}
\end{array}\right.
\end{equation*}

It is easy to see that the two-parameter family $\{L_{s,t}\}$
is a continuous family of bounded selfadjoint operators
on $L_2(M,\mathbb{E})$.

We define $\{\cal{A}_{s,t}\}$ a two-parameter family of operators on
$H^1(M,\mathbb{E})$ defined as
\begin{equation*}
\cal{A}_{s,t}(f)=\left\{\begin{array}{ll}
(A+C_s)(f) & \text{on $M_{-}$},\\
(A+C_t)(f) & \text{on $M_{+}$},
\end{array}\right.
\end{equation*}
then since $\cal{A}_{s,t}=A+C_s+L_{s,t}$ we have
\begin{prop}
Each operator $\cal{A}_{s,t}$ satisfies the `a priori' estimate
\begin{equation*}
\norm{u}_1\le C(\norm{\cal{A}_{s,t}(u)}_0+\norm{u}_0)
\end{equation*}
with a uniform constant $C>0$.

Here we denote $\norm{\cdot}_1$ the first order Sobolev norm of
$u\in H^1(M,\mathbb{E})$.
\end{prop}

By this proposition we see that the direct sum $\Lambda_s^{-}\oplus
\Lambda_t^{+}$ of the Cauchy data spaces of the operator
$(\cal{A}_0^{-})^*+C_s$ and $(\cal{A}_0^{+})^*+C_t$ 
give rise to a Fredholm pair in 
$\boldsymbol\beta^{-}\oplus \boldsymbol\beta^{+}$ and the family 
$\{\Lambda_s^{-}\oplus\Lambda_t^{+}\}$ 
is a continuous two-parameter family of Lagrangian
subspaces in $\boldsymbol\beta^{-}\oplus \boldsymbol\beta^{+}$.  Hence we have

\begin{thm}\label{thm:pre-main}
$\Sf\{A+C_t\}=\Sf(\{\cal{A}_{s,0}\}_{s\in I})+\Sf(\{\cal{A}_{1,t}\}_{t\in I})$.
\end{thm}

This formula can be seen as a splitting formula of 
the spectral flow for the family $\{A+C_t\}$. 
However each term
in the right hand side is a spectral flow corresponding to the
operator family $\{\cal{A}_{s,0}\}_{s\in I}$ and $\{\cal{A}_{1,t}\}_{t\in I}$
on the whole manifold $M$.
The purpose of this paper is to study when we can express these spectral flows
in terms of spectral flows of families of operators restricted to 
the separated submanifolds $M_{\pm}$. 
We are not sure whether such a 
kind of splitting formula holds for any family 
of selfadjoint elliptic operators. 
One difficulty for expressing each
term as a spectral flow of a family defined on the component
$M_{\pm}$ lies  in the general
spectral flow formula.  
Because there we have a formula for the spectral flow
in terms of the Maslov index in the distribution space $\boldsymbol\beta$ 
and these two spaces of boundary values $\boldsymbol\beta^{\pm}$ taken
from each side of the separating hypersurface $\Sigma$ are not
identical even for the case of the operators 
of the product form in a neighborhood of $\Sigma$.  
So it will be natural to expect that 
if we have a reformulation of the  spectral flow formula 
in terms of the $L_2$-space on the hypersurface $\Sigma$, 
then we will have such a splitting formula.
We would like to make these situations clearer.
So we recall in the next subsection 
a reduction formula for the Maslov index in the infinite dimension.

\subsection{A reduction formula for the Maslov index}

Let $\boldsymbol\beta$ and $L$ be two symplectic Hilbert spaces with the symplectic
form $\omega_{\boldsymbol\beta}$ and $\omega_L$.  We assume that each space
$\boldsymbol\beta$ and $L$ decomposes as a direct sum of two Lagrangian subspaces 
in the following way:
\begin{equation*}
\boldsymbol\beta=\theta_{-}\oplus \theta_{+}
\end{equation*}
and
\begin{equation*}
L=L_{-}\oplus L_{+}.
\end{equation*}
We also assume that there are injective maps
\begin{align*}
\mathbf{i}_{-}&:\; \theta_{-}\rightarrow L_{-}\\
\mathbf{i}_{+}&:\; L_{+}\rightarrow \theta_{+}
\end{align*}
such that 
$\omega_L(\mathbf{i}_{-}(x),y)=\omega_{\boldsymbol\beta}(x,\mathbf{i}_{+}(y))$
for $x\in \theta_{-}$ and $y\in L_{+}$, and moreover we assume that the
images under these maps are dense.

Then we have

\begin{thm}[\cite{BFO}]\label{thm:reduction}
There is a natural continuous map
\begin{equation*}
\ell:\;\cal{FL}_{\theta_{-}}(\boldsymbol\beta)\rightarrow\cal{FL}_{L_{-}}(L)
\end{equation*}
satisfying
\begin{equation*}
\Mas(\{\mu_t\},\theta_{-})=\Mas(\{\ell(\mu_{t})\},L_{-}).
\end{equation*}
\end{thm}

\subsection{Splitting formula for a spectral flow}

{}From this section on we assume that the operators $A+C_s$
appearing in the family are all of the product form near the
hypersurface $\Sigma$, that is, the
operator $A+C_s$ has the following form on a cylindrical
neighborhood $U(\Sigma)\cong (-1,1)\times\Sigma$ of $\Sigma$:     
\begin{equation*}
A+C_s=\sigma\left(\frac{\partial}{\partial \tau}+B_s\right)
\end{equation*}
where $\tau\in(-1,1)$, $\sigma$ is a bundle map of
$\mathbb{E}|_{\Sigma}$ defined on $\Sigma$, $B_s$ is 
a selfadjoint elliptic operator on 
$\mathbb{E}|_{\Sigma}$,
$\pi^{*}(\mathbb{E}|_{\Sigma})\cong \mathbb{E}|_{U(\Sigma)}$, and
$\pi:\;U(\Sigma)\rightarrow\Sigma$ is defined through a fixed identification 
$(-1,1)\times\Sigma\cong U(\Sigma)$.

Let $\{\varphi_k\}_{k\in\Z\setminus\{0\}}$ be a complete orthonormal
system of eigen-sections of the operator $B_0$, then we have
a concrete characterization of the boundary value space
$\boldsymbol\beta^{\pm}$ as follows.  

\begin{prop}[\cite{BF2}]\label{prop:4.6}
\begin{align*}
\boldsymbol\beta^{+}&=
{\overline{[\{\varphi_k\}]}}_{k<0}^{H^{1/2}(\Sigma)}
\oplus{\overline{[\{\varphi_k\}]}}_{k>0}^{H^{-1/2}(\Sigma)}
=\cal{B}^+_{-}\,+\,\cal{B}^+_{+}\\
\boldsymbol\beta^{-}&=
{\overline{[\{\varphi_k\}]}}_{k<0}^{H^{-1/2}(\Sigma)}
\oplus{\overline{[\{\varphi_k\}]}}_{k>0}^{H^{1/2}(\Sigma)}
=\cal{B}^-_{-}\,+\,\cal{B}^-_{+}\\
\end{align*}
where we mean by $\cal{B}^{+}_{-}$ = 
${\overline{[\{\varphi_k\}]}}_{k<0}^{H^{1/2}(\Sigma)}$
the completion of the space $[\{\varphi_k\}_{k<0}]$
spanned by $\{\varphi_k\}_{k<0}$ with respect to the Sobolev norm
$H^{1/2}(\Sigma)= H^{1/2}(\Sigma, \mathbb{E}|_{\Sigma})$ and so on.
\end{prop}

Let us consider the family of operators $\{\cal{A}_{s,0}\}$
\begin{equation}\label{equation:two-para}
\cal{A}_{s,0}(f)=\left\{\begin{array}{ll}
(A+C_s)(f) & \text{on $M_{-}$},\\
(A+C_0)(f) & \text{on $M_{+}$},
\end{array}\right.
\end{equation}
for $f\in H^1(M,\mathbb{E})$.
The spectral flow of this family is equal to
$\Mas(\{\Lambda_s^{-}\oplus\Lambda_0^{+}\},\boldsymbol\delta)$
by our general spectral flow formula (Theorem \ref{thm:general}).

Then let us take $\boldsymbol\beta$ and $L$ in the Theorem
\ref{thm:reduction} as 
\begin{align*}
\boldsymbol\beta &= \boldsymbol\beta^{-}\oplus\boldsymbol\beta^{+}\\
L &=L_2(\Sigma,\mathbb{E}|_{\Sigma})
\oplus 
L_2(-\Sigma,\mathbb{E}|_{\Sigma})
=L_2(\Sigma,\mathbb{E}|_{\Sigma})\boxplus
L_2(\Sigma,\mathbb{E}|_{\Sigma}).
\end{align*}
Here 
$L_2(-\Sigma) = L_2(-\Sigma,\mathbb{E}|_{\Sigma})$ 
means the symplectic Hilbert space with the
opposite sign of symplectic structure of $L_2(\Sigma,\mathbb{E}|_{\Sigma})$, or the
orientation of $\Sigma$ should be reversed.

We have decompositions of $\boldsymbol\beta=\theta_{-}\oplus\theta_{+}$ and
$L=L_{-}\oplus L_{+}$ as follows:
\begin{align*}
\theta_{-}& = \boldsymbol\delta =\{(\varphi,\varphi)\in 
\boldsymbol\beta^{-}\oplus \boldsymbol\beta^{+}\;|\;{}^{\exists}
f\in H^1(M,\mathbb{E}),\gamma_{\pm}(f|_{M_{\pm}})=\varphi\},\\
\theta_{+}& =\{(x,0,0,y)\,\in\,
\cal{B}^-_{-}\,+\,\cal{B}^-_{+}\,+\,\cal{B}^+_{-}\,+\,\cal{B}^+_{+}\,|\,x\in
\overline{[\{\varphi_k\}]}_{k<0}^{H^{-1/2}(\Sigma)},\,
y\in \overline{[\{\varphi_k\}]}_{k>0}^{H^{-1/2}(\Sigma)}\},\\
L_{-}&=\Delta=
\text{the diagonal of}\,\, 
L_2(\Sigma,\mathbb{E}|_{\Sigma})\boxplus L_2(\Sigma,\mathbb{E}|_{\Sigma}),\\
L_{+}&=\{(u,0,0,v)|\,u\in\overline{[\{\varphi_k\}]}_{k<0}^{L_2(\Sigma)},\,
v\in \overline{[\{\varphi_k\}]}_{k>0}^{L_2(\Sigma)}\}.
\end{align*}

We have obvious embeddings 
$\mathbf{i}_{-}\, :\, \theta_{-}\rightarrow L_{-}$
and
$\mathbf{i}_{+}\,:\, L_{+}\rightarrow \theta_{+}$
in the space of distributions on $\Sigma$ 
satisfying the properties as required in Theorem \ref{thm:reduction}.
Then we have a spectral flow formula in the $L_2$ setting, since
in this case the map $\ell$ means just taking the 
intersection in the distribution space on $\Sigma$:

\begin{prop}
$\Sf\{\cal{A}_{s,0}\}=\Mas(\{(\Lambda_s^{-}\cap L_2(\Sigma))
\boxplus(\Lambda_0^{+}\cap L_2(\Sigma))\},\Delta)$.
\end{prop}

Now from Proposition \ref{prop:pair} we have
\begin{prop}
\begin{align*}
&\Mas(\{(\Lambda_s^{-}\cap L_2(\Sigma))
\boxplus(\Lambda_0^{+}\cap L_2(\Sigma))\},\Delta)\\
&=\Mas(\{\Lambda_s^{-}\cap L_2(\Sigma)\},\Lambda_0^{+}\cap L_2(\Sigma)).
\end{align*}
\end{prop}

Next let us consider the family of operators $\{\cal{T}_s^{-}\}$ on 
$M_{-}$ with the domain $D_0$,
\begin{equation}\label{def:leftelliptic}
D_0=\{f\in H^1(M_{-},\mathbb{E}|_{M_{-}})\;|\;{}^{\exists}\varphi\in\Lambda_0^{+},
\gamma_{-}(f)=\varphi\}
\end{equation}
and
\begin{equation*}
\cal{T}_s^{-}(f)=(A+C_s)(f)\quad\text{on $M_{-}$}.
\end{equation*}

{}From the reduction theorem \ref{thm:reduction}, 
$\gamma_{-}(D_0)$ is a Lagrangian subspace
in $\boldsymbol\beta^{-}$. Hence we have
\begin{prop}
Each operator $\cal{T}_s^{-}$ on $D_0$ is a selfadjoint operator. So this
implies that 
\begin{equation*}
\norm{u}_1\le c'(\norm{\cal{T}_s^{-}u}_0+\norm{u})\quad\text{for $u\in D_0$}
\end{equation*}
with a uniform constant $c'>0$. This inequality implies that we have
well-defined spectral flow for this continuous family
$\{\cal{T}_s^{-}\circ\sqrt{1+(\cal{T}_s^{-})^2}^{-1}\}$.
\end{prop}

Hence we have
\begin{prop}
\begin{align*}
\Sf(\{\cal{T}_s^{-}\})&=\Mas(\{\Lambda_s^{-}\},\gamma_{-}(D_0))\\
&=\Mas(\{\Lambda_s^{-}\cap L_2(\Sigma)\},\Lambda_0^{+}\cap L_2(\Sigma))\\
&=\Sf\{\cal{A}_{s,0}\}.
\end{align*}
\end{prop}

Similarly when we define the family of operators $\{\cal{T}_t^{+}\}$
on $M_{+}$ with the domain $D_1$ given by
\begin{equation}\label{def:rightelliptic}
D_1=\{v\in H^1(M_{+},\mathbb{E}|_{M_+})\;|\;{}^{\exists}\varphi\in\Lambda_1^{-},
\gamma_{+}(v)=\varphi\}
\end{equation}
and
\begin{equation*}
\cal{T}_t^{+}(v)=(A+C_t)v\quad\text{on $M_{+}$}.
\end{equation*}
then we have

\begin{prop}
\begin{align*}
\Sf(\{\cal{T}_t^{+}\})&=\Mas(\{\Lambda_t^{+}\},\gamma_{+}(D_1))\\
&=\Mas(\{\Lambda_t^{+}\cap L_2(\Sigma)\},\Lambda_1^{-}\cap L_2(\Sigma))\\
&=\Sf\{\cal{A}_{0,t}\}.
\end{align*}
\end{prop}

Summing up the results above we have our main Theorem of a splitting
formula for the spectral flow $\Sf\{A+C_t\}$.

\begin{thm}\label{thm:main}
$\Sf\{A+C_t\}=\Sf\{\cal{T}_t^{-}\}+\Sf\{\cal{T}_t^+\}$.
\end{thm}

\begin{rem}
\begin{enumerate}
\item
We have proved Theorem \ref{thm:main} under the assumption of the 
product form structure for the operators near the hypersurface
$\Sigma$. It is important that this assumption implies that 
the conditions 
to define the domains $D_0$ and $D_1$ 
become selfadjoint elliptic boundary conditions.
\item
If the operator is not of product form near 
the hypersurface $\Sigma$, but if an operator $A+C_0$
is invertible, then we can prove that the condition
in (\ref{def:leftelliptic}) gives us an elliptic selfadjoint
boundary condition.

\begin{proof}
The invertibility of the operator $A+C_0$ is equivalent
to the transversality of the two Lagrangian subspaces
$\boldsymbol\delta$ and $\Lambda^{-}_0\oplus\Lambda^{+}_0$ in 
$\boldsymbol\beta = \boldsymbol\beta^{-}\oplus\boldsymbol\beta^{+}$.
Hence any $(z,0)\in \boldsymbol\beta^{-}\oplus\boldsymbol\beta^{+}$
can be written as 
$(z,0) =(x,y) + (a,a)\in\Lambda^{-}_0\oplus\Lambda^{+}_0 
+ \boldsymbol\delta$.  
So $a = -y\in H^{1/2}(\Sigma,\mathbb{E}|_{\Sigma})\cap\Lambda^+_0$ and
$z=x-y$, that is,
$\boldsymbol\beta^- =\Lambda^-_0 + \gamma(D_0)$.  
Hence this implies $\gamma(D_0)$ must be a Lagrangian subspace,
because $\gamma(D_0)$ is always isotropic, and we have also
that the closedness of the operator $\cal{T}^-_0$ 
defined on $D_0 \subset
H^1(M,\mathbb{E})$ gives us the ellipticity of the condition
for any operator $\cal{T}^-_t= (A+C_0)+(C_t-C_0)$, $t\in I$.
\end{proof}
\end{enumerate}
\end{rem}

By the above remark, if the operator $A+C_0$
is invertible we have
\begin{prop}
$\Sf\{\cal{T}_{s}\} = \Mas(\{\Lambda^-_s\},\gamma(D_0))$.
\end{prop}

However we do not know whether $\Mas(\{\Lambda^-_s\},\gamma(D_0))$
coincides with
$\Mas(\{\Lambda^-_s\oplus\Lambda^+_0\},\boldsymbol\delta)$
without the assumption of the product form near $\Sigma$.


\appendix

\setcounter{secnumdepth}{1}

\setcounter{thm}{0}
\setcounter{section}{1}
\setcounter{equation}{0}

\section*{Appendix}

In this appendix we explain the space $\cal{U}_{\cal{F}}(H_J)$ in
the framework of the complexified symplectic Hilbert space
(Proposition \ref{prop:introduce})
and give a proof for the isomorphisms:
\begin{prop}\label{prop:homotopy}
$\pi_{1} (\mathcal{U}_{\lambda}(H_J))
\underset{(\cal{S}_{\lambda})_*}{\xrightarrow{\sim}}\pi_1(\cal{U}_{\cal{F}}(H_J))
\xrightarrow{\sim}\Z$.
\end{prop}

Let $H$ be a separable symplectic Hilbert space with the
symplectic form $\omega$, an inner product $\inner{\cdot}{\cdot}$
and an almost complex structure  $J$.  
We assume that 
all of these three are compatible as before, that is,
\begin{equation*}
\omega(x,y) =\inner{J(x)}{y},\quad J^2 = -Id.
\end{equation*}
Let $\mathcal{L}^{\C}(H\otimes \C)$ be the space of complex Lagrangian
subspace in $H\otimes \C$:
\begin{equation*}
\cal{L}^{\C}(H\otimes\C)=\{l\;|\;
\text{$l$ is a subspace such that $l^{\perp}=J(l)$}\}.
\end{equation*}
Then the subgroup of the unitary operators in $H\otimes \C$,
denoted by
$\mathcal{U}_0(H\otimes \C)$, 
consisting of those operators $U$ such that 
$U(l)^{\perp}=J(U(l))$ for any $l$ in $\mathcal{L}^{\C}(H\otimes \C)$
acts on $\mathcal{L}^{\C}(H\otimes \C)$ transitively.
This condition for $U \in\mathcal{U}_0(H\otimes \C)$ is
equivalent
to the statement that it commutes with the almost
complex structure $J$ (notice that 
$J$ should be considered to be complexified).

Taking the complexification of $\lambda\in \cal{L}(H)$ gives us a
natural embedding
$\cal{L}(H)\rightarrow\cal{L}^{\C}(H\otimes\C)$,
and its restriction to $\mathcal{FL}_{\lambda}(H)$ has the
image in $\mathcal{FL}_{\lambda\otimes \C}^{\C}(H\otimes \C)$, a
subspace of 
$\cal{L}^{\C}(H\otimes\C)$ consisting of those subspaces which are Fredholm
pairs with $\lambda \otimes \C$.
We denote this map by $\cal{C}$.

When we consider an operator $U\in \mathcal{U}(H_J)$ as a real operator and take
its
complexification, we denote it by $U^{\C}$, then
$U^{\C}$ is in $\mathcal{U}_0(H\otimes \C)$ and 
we have $U(\mu)\otimes\C$ = ${U^{\C}}(\mu\otimes\C)$, 
$\mu\in \cal{L}(H)$.

Let $E_{\pm}=\{z\in H\otimes\C\;|\; J(z)=\pm \sqrt{-1}z\}$, then we have
an orthogonal decomposition of $H\otimes\C$ as
\begin{equation*}
H\otimes\C=E_{+}\oplus E_{-},
\end{equation*}
and if $U\in\cal{U}_0(H\otimes\C)$, then
$U(E_{\pm})=E_{\pm}$.
Hence we have an isomorphism
\begin{equation*}
\cal{U}_0(H\otimes\C)\cong \cal{U}(E_{+})\times \cal{U}(E_{-}),
\end{equation*}
where $\cal{U}(E_{+})$ denotes the group of unitary operators on $E_{+}$, and
so on.
Also the space $\cal{L}^{\C}(H\otimes\C)$ is identified
with the space of graphs  of unitary transformation 
$U$ $\in$ $\cal{U}(E_{+},E_{-}),\,\, U:E_{+}\rightarrow E_{-}$. 

Let $\frak{K}: H_J \rightarrow E_{+}$, $u \mapsto u\otimes 1
- J(u)\otimes\sqrt{-1}$ and 
$\frak{k}:H_J \rightarrow E_{-}$, 
$ u \mapsto u\otimes 1 + J(u)\otimes\sqrt{-1}$,
be an isomorphism and an anti-isomorphism 
then
\begin{lem}
The following diagram is commutative.
\begin{equation*}
\begin{CD}
H_J @>{\tau_{\lambda}}>>H_J\\
@V{\frak{K}}VV  @VV{\frak{k}}V\\
E_{+} @>>{T_{\lambda}}>E_{-},\\
\end{CD}
\end{equation*}
where $\tau_{\lambda}$ is the complex conjugation defined
through the identification $H_J \cong \lambda\otimes\C$, and
the graph of the unitary operator $T_{\lambda}$ is 
$\lambda\otimes\C = \{x+T_{\lambda}(x)\,|\,x\in E_{+}\}$.
\end{lem}

Now we have 
\begin{prop}\label{prop:introduce}
Let 
$\varPhi : \cal{U}_{\cal{F}}(H_J)\rightarrow 
\cal{FL}_{\lambda\otimes\C}^{\C}(H\otimes\C)$
be a map defined by
$\varPhi(V) = \text{the graph of the unitary operator}\, 
\, -\frak{k}\circ V \circ \tau_{\lambda}\circ
\frak{K}^{-1} \in \cal{U}(E_+,E_{-})$.
Then 
$\varPhi$ is an isomorphism 
and
the following diagram is commutative:
\begin{equation*}
\begin{psmatrix}
[name=FL] \cal{FL}_{\lambda}(H)  & [name=FLC]
\cal{FL}_{\lambda\otimes\C}^{\C}(H\otimes\C)\\
 & [name=U] \cal{U}_{\cal{F}}(H_J)
\psset{arrows=->}
\ncline{FL}{FLC}^{\cal{C}}\ncline{FL}{U}_{\cal{S}_{\lambda}}
\ncline{U}{FLC}>{\varPhi}.
\end{psmatrix}
\end{equation*}
\end{prop}

\begin{proof}
We only prove the commutativity of the diagram.
Let $U \in \cal{U}_{\lambda}(H_J)$.  Since $U^{\C}|_{E_{\pm}}$ can be
identified with $U$ through the map $\frak{K}$ and $\frak{k}$
respectively, 
we have  
$U^{\C}(\lambda^{\perp}\otimes\C)$ = $\{U(x)-U\circ
T_{\lambda}(x)\,|\,x\in E_{+}\}$ = $\{x-U\circ T_{\lambda}\circ
U^{-1}(x)\,|\, x\in E_{+}\}$.  
By the above lemma
$\frak{k}\circ U\circ T_{\lambda}\circ U^{-1}\circ \frak{K}^{-1}$
= $\frak{k}\circ U\circ \tau_{\lambda}\circ U^{-1}\circ
\tau_{\lambda}\circ \tau_{\lambda}\circ \frak{K}^{-1}$ =
$\frak{k}\circ U\circ \theta_{\lambda}(U) \circ\tau_{\lambda}\circ \frak{K}^{-1}$,
which gives the commutativity of the diagram.
\end{proof}

Let $W$ be a closed finite codimensional subspace in 
$\lambda\otimes\C$ and we denote by $\cal{FL}_W^{(0)}$ 
a subspace of $\cal{FL}_{\lambda\otimes\C}^{\C}(H\otimes\C)$
consisting of those subspaces $l$ which do not intersect with $W$.
Let $H_W$ =
$J(W^{\perp}\cap(\lambda\otimes\C))+W^{\perp}\cap(\lambda\otimes\C)$,
and $\cal{L}(H_W)$ be the similar space as $\cal{L}(H\otimes\C)$
(note $H_W$ is invariant under the map $J$).  
$\cal{L}(H_W)$ is identified with the space of unitary operators on 
$W^{\perp}\cap(\lambda\otimes\C)$.  
Let 
${\pi}_W : 
\cal{FL}_W^{(0)} \ni l \rightarrow 
(l\,\cap\,(J(W^{\perp}\cap(\lambda\otimes\C))+\lambda\otimes\C)+W)\cap
W^{\perp} \in \cal{L}(H_W)$, 
and then ${\pi}_W : \cal{FL}^{(0)}_W \rightarrow \cal{L}(H_W)$ is a fiber bundle 
with contractible fibers.   
A typical fiber 
= $\pi^{-1}_W(J((\lambda\otimes\C)\cap W^{\perp}))$ is isomorphic to
the space $\widehat{\cal{B}}(W) \times \cal{B}(W, (\lambda\otimes\C)\cap
W^{\perp})$, where $\widehat{\cal{B}}(W)$ is the space of selfadjoint
operators on $W$ and $\cal{B}(W, (\lambda\otimes\C)\cap W^{\perp})$ is the space
of bounded operators from $W$ to $(\lambda\otimes\C)\cap W^{\perp}$.
Unfortunately  
for any pair of such subspaces $W_1$ and $W_2$ satisfying $W_1\subset W_2$
there are no natural map $\cal{L}(H_{W_2}) \rightarrow
\cal{L}(H_{W_1})$ which
makes the diagram
\begin{equation*}
\begin{CD}
\cal{FL}_{W_2}^{(0)} @>>> \cal{FL}_{W_1}^{(0)}\\
@VV{\pi_{W_2}}V @VV{\pi_{W_1}}V\\
\cal{L}(H_{W_2}) @>>> \cal{L}(H_{W_1}) 
\end{CD}
\end{equation*}
commutative. However if we define a map $s_W : \cal{L}(H_W)
\rightarrow \cal{FL}_W^{(0)}(H\otimes\C)$ 
by $s_W(l) = l+J(W)$, then $\pi_W\circ s_W =Id$ and 
we have the following commutative diagram:
\begin{equation*}
\begin{CD}
\cal{L}(H_{W_2}) @>>{\bf{i}}_{W_1,W_2}> \cal{L}(H_{W_1}) \\
@VV{s_{W_2}}V @VV{s_{W_1}}V\\
\cal{FL}_{W_2}^{(0)} @>>> \cal{FL}_{W_1}^{(0)},
\end{CD}
\end{equation*}
where the map ${\bf{i}}_{W_1,W_2}:\cal{L}(H_{W_2}) \rightarrow \cal{L}(H_{W_1})$
is defined as ${\bf{i}}_{W_1,W_2}(l) = l+J(W_2 \cap W_1^{\perp})$.

Then 
for any compact subset $K$ in
$\cal{FL}_{\lambda\otimes\C}^{\C}(H\otimes\C)$
we can find such a finite codimensional subspace $W$ in $\lambda\otimes\C$ 
that for any $l$ in $K$, $l\cap W = \{0\}$, so  
$\bigcup\cal{FL}_W^{(0)}$ 
= $\cal{FL}_{\lambda\otimes\C}^{\C}(H\otimes\C)$.
Hence 
$\lim\limits_{W\rightarrow \{0\}} \pi_{k}(\cal{FL}_W^{(0)})$ =
$\pi_k(\cal{U}_{\cal{F}}(H_J))$. 
These facts show that the homotopy groups of
$\cal{U}_{\cal{F}}(H_J)$ coincide with the stable homotopy groups
of unitary groups, which together gives the proof of
Proposition \ref{prop:homotopy}.




\end{document}